\documentclass[11pt]{article}
\usepackage{geometry}
\usepackage{amssymb, amsmath, amsthm, fullpage, color, tikz}
\newtheorem{theorem}{Theorem}[section]

\newtheorem{definition}[theorem]{Definition}

\numberwithin{equation}{section}

\begin{document}

\title{Parking Cars of Different Sizes}
\markright{Parking Cars of Different Sizes}
\author{Richard Ehrenborg and Alex Happ}
\date{}

\maketitle

\begin{abstract}
We extend the notion of parking functions
to parking sequences, which include
cars of different sizes, and prove
a product formula for the number of such sequences.
\end{abstract}

\section{The result.}

Parking functions were
first
introduced by
Konheim and Weiss~\cite{Konheim_Weiss}.
The original concept was that of a linear parking lot 
with $n$ available spaces, and 
$n$~cars with a stated parking preference. 
Each car would, in order, 
attempt to park in its preferred spot. 
If the car found its preferred spot occupied, 
it would move to the next available slot. 
A parking function is a sequence 
of parking preferences that would allow 
all $n$ cars to park according to this rule. 
This definition is equivalent 
to the following formal definition:
\begin{definition}
   Let $\vec{a}=(a_1,a_2,\dots,a_n)$
   be a sequence of positive integers,
   and let $b_1\leq b_2\leq\cdots\leq b_n$
   be the increasing rearrangement of $\vec{a}$.
   Then the sequence $\vec{a}$ is a parking function
   if and only if
   $b_i\leq i$ for all indexes $i$.
\end{definition}

It is well known that the number of 
such parking functions is $(n+1)^{n-1}$.
This is Cayley's formula for the number
of labeled trees on $n+1$ nodes
and Foata and Riordan found
a bijective proof~\cite{Foata_Riordan}. 
Stanley discovered the relationship 
between parking functions 
and non-crossing partitions~\cite{Stanley_I}. 
Further connections have been found to other structures, 
such as
priority queues~\cite{Gilbey_Kalikow},
Gon\v{c}arov polynomials~\cite{Kung_Yan_I}
and
hyperplane arrangements~\cite{Stanley_II}.

The notion of a parking function
has been generalized in myriad ways; see the sequence of
papers~\cite{Chebikin_Postnikov,Kung_Yan_I,Kung_Yan_II,Kung_Yan_III,Yan}.
We present here a different generalization,
returning to the original idea of parking cars.
This time the cars have different sizes, and
each takes up a number of adjacent parking spaces.
\begin{definition}
Let there be $n$ cars $C_{1},\dots,C_n$
of sizes $y_{1},\dots,y_{n}$, 
where $y_{1}, \ldots, y_{n}$ are positive integers.
Assume there are $\sum_{i=1}^{n} y_{i}$ spaces in a row.
Furthermore, let car $C_{i}$ have the preferred spot $c_{i}$.
Now let the cars
in the order $C_{1}$ through $C_{n}$
park according to the following rule:
\begin{quote}
Starting at position $c_{i}$,
car $C_{i}$
looks for the first empty spot $j \geq c_{i}$.
If the spaces~$j$ through $j+y_{i}-1$ are empty, then car $C_{i}$
parks in these spots.
If any of the spots
$j+1$ through $j+y_{i}-1$ is already occupied,
then there will be a collision, and the result is not a parking sequence.
\end{quote}
Iterate this rule for all the cars
$C_{1}, C_{2}, \ldots, C_{n}$.
We call $(c_{1},\dots, c_n)$ a \emph{parking sequence}
for $\vec{y}=(y_{1},\dots,y_{n})$
if all $n$ cars can park without any collisions
and without leaving
the $\sum_{i=1}^{n} y_{i}$ parking spaces.
\end{definition}

As an example, consider three cars of sizes
$\vec{y}=(2,2,1)$
with preferences $\vec{c}=(2,3,1)$.
Then there are $2+2+1=5$ available parking spaces,
and the final configuration of the cars is
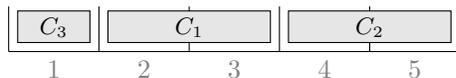
\[\begin{tikzpicture}[scale=1.2]
\draw (0,0) -- (5,0);
\foreach \x in {0,...,5}
\draw (\x,0) -- (\x,0.5);
\foreach \x in {1,...,5}
\node[gray] at (\x-0.5,-0.2) {\small$\x$};
\draw[fill=gray!20] (0.1,0.1) rectangle (0.9,0.45);
\draw[fill=gray!20] (1.1,0.1) rectangle (2.9,0.45);
\draw[fill=gray!20] (3.1,0.1) rectangle (4.9,0.45);
\node at (0.5,0.265) {\footnotesize $C_{3}$};
\node at (2,0.265) {\footnotesize $C_{1}$};
\node at (4,0.265) {\footnotesize $C_{2}$};
\end{tikzpicture}\]
All cars are able to park, so this yields a parking sequence.

There are two ways in which a sequence can fail
to be a parking sequence.
Either a collision occurs, or a car passes the end of the parking lot.
As an example, consider three cars with
$\vec{y}=(2,2,2)$
and preferences $\vec{c}=(3,2,1)$.
Then we have $2+2+2=6$ parking spots,
and the first car parks in its desired spot:
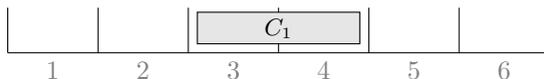
\[\begin{tikzpicture}[scale=1.2]
\draw (0,0) -- (6,0);
\foreach \x in {0,...,6}
\draw (\x,0) -- (\x,0.5);
\foreach \x in {1,...,6}
\node[gray] at (\x-0.5,-0.2) {\small$\x$};
\draw[fill=gray!20] (2.1,0.1) rectangle (3.9,0.45);
\node at (3,0.265) {\footnotesize $C_{1}$};
\end{tikzpicture}\]
However, the second car prefers spot $2$,
and since spot $2$ is open, he tries to take spots~$2$ and~$3$,
but collides with $C_{1}$ in the process. Hence, this is not a
parking sequence.

If, instead, we had $\vec{y}=(2,2,2)$
and
$\vec{c}=(2,5,5)$,
then again the first two cars are able to park with no difficulty:
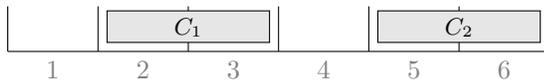
\[\begin{tikzpicture}[scale=1.2]
\draw (0,0) -- (6,0);
\foreach \x in {0,...,6}
\draw (\x,0) -- (\x,0.5);
\foreach \x in {1,...,6}
\node[gray] at (\x-0.5,-0.2) {\small$\x$};

\draw[fill=gray!20] (1.1,0.1) rectangle (2.9,0.45);
\draw[fill=gray!20] (4.1,0.1) rectangle (5.9,0.45);

\node at (2,0.265) {\footnotesize $C_{1}$};
\node at (5,0.265) {\footnotesize $C_{2}$};

\end{tikzpicture}\]
But car $C_{3}$ will pass by all the parking spots
after his preferred spot without seeing an empty spot.
Hence, this also fails to be a parking sequence.

The classical notion of
parking function is obtained when
all the cars have size~$1$, that is,
$\vec{y}=(1,1, \ldots, 1)$.
Note in this case that there are no possible collisions.

In the classical case, any permutation of a
parking function is again a parking function.
This is not true for cars of larger size.
As an example, note for $\vec{y} = (2,2)$
that $\vec{c} = (1,2)$ is a parking sequence.
However, the rearrangement
$\vec{c}\,^{\prime} = (2,1)$ is not a parking sequence.
This shows that the notion of parking sequence
differs from the notion of
parking function in the
papers~\cite{Chebikin_Postnikov,Kung_Yan_I,Kung_Yan_II,Kung_Yan_III,Yan}.

The classical result is that the
number of parking functions is given
by $(n+1)^{n-1}$;
see~\cite{Konheim_Weiss}.
For cars of bigger sizes we have the following result:
\begin{theorem}
The number of parking sequences $f(\vec{y})$ for car sizes
$\vec{y}=(y_{1},\dots,y_n)$
is
given by the product
\[
f(\vec{y})=
(y_{1}+n)
\cdot (y_{1}+y_{2}+n-1)
\cdots (y_{1}+\cdots+y_{n-1}+2).\]
\label{theorem_parking}
\end{theorem}

\section{Circular parking arrangements.}

Consider $M = y_{1} + y_{2} + \cdots + y_{n}  + 1$ parking spaces
arranged in a circle.
We will consider parking cars on this circular arrangement,
without a cliff for cars to fall off.
Observe that when all the cars
have parked, there will be one empty spot left over.
We claim that there are
\begin{equation}
M \cdot f(\vec{y})
=
(y_{1}+n)
\cdot (y_{1}+y_{2}+n-1)
\cdots 
(y_{1}+\cdots+y_{n}+1) .
\label{equation_circular}
\end{equation}
such circular parking sequences.
The first car $C_{1}$ has $M$
ways to choose its parking spot.

The next step is counterintuitive.
After car $C_{1}$ has parked, erase the markings for the remaining $y_{2}+ \cdots + y_{n} + 1$ spots
and put in $n+1$ dividers. These dividers create $n+1$ intervals on the circle,
where one interval is taken up by $C_{1}$.
Furthermore, these dividers are on wheels and can freely move along the circle.
Each interval will accept one (and only one) car.
For example, consider the case where $n=5$ and $\vec{y}=(2,5,1,3,2)$ so that $M=2+5+1+2+3+1=14$, and $c_1=5$.

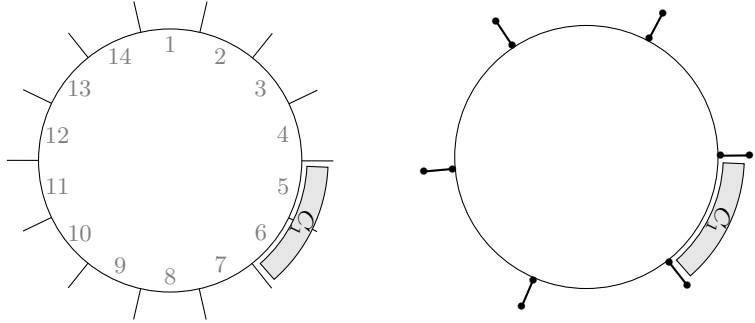
\[\begin{tikzpicture}[scale=0.7]
	\draw (0,0) circle (2.5);
	\foreach \x in {0,...,13}
		\draw (102.85+\x*25.7:2.5) -- (102.85+\x*25.7:3.1);
	\foreach \x in {1,...,14}
		\node at (115.7-\x*25.7:2.2) {\color{gray}\footnotesize $\x$};
	\draw[fill=gray!20] (-2.45:3) arc (-2.45:-48.85:3) -- (-48.85:2.6) arc (-48.85:-2.45:2.6) -- cycle;
		\node[rotate=-115.65] at (-25.65:2.795) {\footnotesize $C_1$};
  \end{tikzpicture}
  \hspace{1cm}
  \begin{tikzpicture}[scale=0.7]
	\draw (0,0) circle (2.5);
	\draw[thick] (0.55:2.55) -- (0.55:3.1);
		\node at (0.55:2.55) {\tiny$\bullet$};
		\node at (0.55:3.1) {\tiny$\bullet$};
	\draw[thick] (-51.85:2.55) -- (-51.85:3.1);
		\node at (-51.85:2.55) {\tiny$\bullet$};
		\node at (-51.85:3.1) {\tiny$\bullet$};
	\foreach \x in {1,...,4}
		\draw[thick] (-51.85-\x*61.52:2.55) -- (-51.85-\x*61.52:3.1);
	\foreach \x in {1,...,4}
		\node at (-51.85-\x*61.52:2.55) {\tiny$\bullet$};
	\foreach \x in {1,...,4}
		\node at (-51.85-\x*61.52:3.1) {\tiny$\bullet$};
		
	\draw[fill=gray!20] (-2.45:3) arc (-2.45:-48.85:3) -- (-48.85:2.6) arc (-48.85:-2.45:2.6) -- cycle;
		\node[rotate=-115.65] at (-25.65:2.795) {\footnotesize $C_1$};
  \end{tikzpicture}\]

We will now create a circular parking sequence,
but only at the end do we obtain the exact positions of cars $C_{2}$ through $C_{n+1}$.
That is, instead of focusing on the number of specific spot preferences each car could have,
we keep track of the order the cars park in,
which will then determine the exact locations of the cars.

The second car has two options. The first is that it has a desired position already taken by $C_{1}$.
In this case, it will cruise until the next empty spot. This can happen in $y_{1}$ ways,
and then car $C_{2}$ obtains the next open interval after the interval $C_{1}$ is in.
Otherwise, the car $C_{2}$ has a preferred spot not already taken. In this case $C_{2}$
has $n$ open intervals to choose from.
The total number of options for $C_2$ is $y_{1} + n$.

The third car $C_{3}$ has the same options. First, it may desire a spot that is already taken,
in which case it will have to cruise
until the next open interval. This can happen in $y_{1} + y_{2}$ ways. 
Note that this count applies to both the case when $C_{1}$ and $C_{2}$ are parked next to each
other, and when $C_{1}$ and $C_{2}$ have open intervals between them. 
Otherwise, $C_{3}$ has $n-1$ open intervals to pick from.

In general, car $C_{i}$ has $y_{1} + \cdots + y_{i-1} + n+2-i$ choices.
This pattern continues up to $C_{n}$, which has
$y_{1} + \cdots + y_{n-1} + 2$ possibilities. For example, suppose $C_2$ and~$C_3$ in our above example have parked as below:

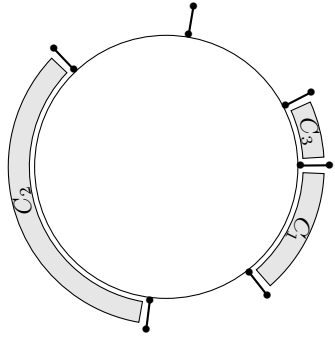
\[\begin{tikzpicture}[scale=0.7]
	\draw (0,0) circle (2.5);
	\draw[thick] (0.55:2.55) -- (0.55:3.1);
		\node at (0.55:2.55) {\tiny$\bullet$};
		\node at (0.55:3.1) {\tiny$\bullet$};
	\draw[thick] (-51.85:2.55) -- (-51.85:3.1);
		\node at (-51.85:2.55) {\tiny$\bullet$};
		\node at (-51.85:3.1) {\tiny$\bullet$};
	\draw[thick] (263:2.55) -- (263:3.1);
		\node at (263:2.55) {\tiny$\bullet$};
		\node at (263:3.1) {\tiny$\bullet$};
	\draw[thick] (133.5:2.55) -- (133.5:3.1);
		\node at (133.5:2.55) {\tiny$\bullet$};
		\node at (133.5:3.1) {\tiny$\bullet$};
	\draw[thick] (27.25:2.55) -- (27.25:3.1);
		\node at (27.25:2.55) {\tiny$\bullet$};
		\node at (27.25:3.1) {\tiny$\bullet$};
	\draw[thick] (80.375:2.55) -- (80.375:3.1);
		\node at (80.375:2.55) {\tiny$\bullet$};
		\node at (80.375:3.1) {\tiny$\bullet$};
		
	\draw[fill=gray!20] (3.55:3) arc (3.55:24.25:3) -- (24.25:2.6) arc (24.25:3.55:2.6) -- cycle;
		\node[rotate=-76.1] at (13.9:2.795) {\footnotesize $C_3$};  
	\draw[fill=gray!20] (-2.45:3) arc (-2.45:-48.85:3) -- (-48.85:2.6) arc (-48.85:-2.45:2.6) -- cycle;
		\node[rotate=-115.65] at (-25.65:2.795) {\footnotesize $C_1$};
	\draw[fill=gray!20] (260:3) arc (260:136.5:3) -- (136.5:2.6) arc (136.5:260:2.6) -- cycle;
		\node[rotate=88.25+15] at (178.25+15:2.795) {\footnotesize $C_2$};
  \end{tikzpicture}\]
Then $C_4$ may either cruise on $C_1$ and $C_3$ (in $y_1+y_3$ ways), 
it may cruise on $C_2$ (in~$y_2$ ways),
or it can pick one of the three available intervals directly. 
In total, $C_4$ has $(y_1+y_3)+y_2+3=11$ ways to park.

One can imagine that when we park a car, we do not set the parking
brake, but put the car in neutral, so that the car
and the dividers can move as necessary to make room for future cars.

Thus the total number of circular parking arrangements of this type is
$$ M
\cdot
  (y_{1} + n)
  \cdot
  (y_{1} + y_{2} + n-1)
  \cdots
  (y_{1}+ \cdots + y_{n-1} + 2)  ,  $$
  where the $i$th factor is the number of options for the car $C_{i}$.
This proves the claim about the number of
circular parking sequences in~\eqref{equation_circular}.

Hence, to prove Theorem~\ref{theorem_parking}
we need only observe that the circular parking
sequences with spot $M$ empty are the same
as our parking sequences.
This follows from the observation that no car
in the circular arrangement has preference $M$,
since otherwise this spot would not be empty.
Furthermore, no car would cruise by this empty spot. 

Observe that the set of circular parking sequences
is invariant under rotation. That is,
if $(c_{1}, c_{2}, \ldots, c_{n})$ is
a parking sequence, then so is
the sequence
$(c_{1}+a, c_{2}+a, \ldots, c_{n}+a)$,
where all the additions are modulo $M$.
In particular, the number of circular parking
sequences with spot $M$ empty
is given by $1/M \cdot M \cdot f(\vec{y}) = f(\vec{y})$.

\section{Concluding remarks.}

The idea of considering a circular arrangement goes
back to Pollak; see~\cite{Riordan}.  
In fact, when all the cars have size $1$,
this argument reduces to his argument
that the number of
classical parking functions is $(n+1)^{n-1}$.

The idea of not using fixed
coordinates when placing
cars in the circular arrangement
is reminiscent of the argument 
Athanasiadis used to compute the characteristic
polynomial of the Shi arrangement~\cite{Athanasiadis}.

\section*{Acknowledgments}

The authors thank two referees for their comments 
as well as Margaret Readdy for her comments
on an earlier draft of this note.
Both authors were partially supported by
National Security Agency grant~H98230-13-1-0280.
The first author wishes to thank the Mathematics Department of
Princeton University where this work was carried out.

\newcommand{\journal}[6]{#1, #2, {\it #3} {\bf #4} (#5), #6.}

\noindent
{\em
Department of Mathematics,
University of Kentucky,
Lexington, KY 40506 \\
{\tt richard.ehrenborg@uky.edu}, {\tt alex.happ@uky.edu}
}

\end{document}